\documentclass[12pt]{article}
\usepackage{amssymb}

\usepackage[latin1]{inputenc}
\def\bea{\begin{eqnarray*} }
\def\eea{\end{eqnarray*} }
\newtheorem{definition}{Definition}

\newtheorem{proposition}[definition]{Proposition}

\newtheorem{theorem}[definition]{Theorem}
\newcommand{\hess}{\mathrm{Hess}}

\newcommand{\rank}{\mathrm{rank}}
\newcommand{\po}{{\hspace*{-1ex}}{\bf .  }}
\newcommand{\R}{\mathbb{R}}
\newcommand{\Les}{\mathbb{L}}

\newcommand{\spa}{\mbox{span}}

\def\e{\epsilon}

\newcommand{\be}{\begin{equation} }
\newcommand{\ee}{\end{equation} }

\def\<{\langle}
\def\>{\rangle}
\def\proof{\noindent{\it Proof: }}
\def\qed{\ifhmode\unskip\nobreak\fi\ifmmode\ifinner\else
\hskip5 pt \fi\fi\hbox{\hskip5 pt \vrule width4 pt
height6 pt  depth1.5 pt \hskip 1pt }}

\begin{document}

\title{Isometric rigidity in codimension two}
\author {Marcos Dajczer \& Pedro Morais}
\date{}
\maketitle


\section{Introduction}

In the local theory of submanifolds a fundamental but difficult problem is to describe  the isometrically deformable isometric immersions  $f\colon\,M^n\to\R^{n+p}$
into Euclidean space  with low codimension $p$ if compared to the dimension $n\ge 3$ of the Riemannian manifold. Moreover, one would like to understand the set of all possible isometric deformations.

 Submanifolds in low codimension  are generically rigid since the fundamental Gauss-Codazzi-Ricci  system of equations  is overdetermined. Rigidity means that there are no other isometric immersion up to rigid motion of the ambient space. As a consequence, it is quite easier to give a generic assumption implying rigidity
than describing the submanifolds that are isometrically deformable.
For instance, the results in \cite{all} and \cite{cd} conclude rigidity provided the second fundamental is ``complicated enough".

 The result stated by Beez \cite{be} in 1876 but correctly proved by Killing \cite{ki} in 1885 says that any deformable hypersurfaces without flat points  has two nonzero principal curvatures (rank two) at any point.  For dimension three the deformation problem for hypersurfaces was first considered by Schur \cite{sc} as early as 1886 and then completely solved  in 1905 by Bianchi \cite{bi}. The general case was solved by
Sbrana \cite{sb} in 1909 and Cartan \cite{ca} in 1916; see \cite{dft}  for additional information. From their
result, we have that even hypersurfaces of rank two are generically rigid.

Outside the hypersurfaces case, the deformation question remains essentially unanswered to this day even for low
codimension  $p=2$. From \cite{cd} or \cite{DF2} any submanifold $f\colon\,M^n\to\R^{n+2}$ is rigid if at any
point the index of relative nullity satisfies $\nu_f\ge n-5$ and any shape operator  has at least three nonzero
principal curvatures. If only the relative nullity condition holds we know from \cite{DG2} that $f$ is genuinely
rigid. This means that given any other isometric immersion $\tilde f\colon\,M^n\to\R^{n+2}$ there is an open
dense  subset of $M^n$ such that restricted to any connected component $f|_U$ and $\hat{f}|_U$ are either
congruent or there are an isometric embedding \mbox{$j\colon\, U\hookrightarrow N^{n+1}$} into a Riemannian
manifold $N^{n+1}$ and either flat or isometric Sbrana-Cartan hypersurfaces $F,\hat F\colon\,N^{n+1}\to\R^{n+2}$
such that  $f|_U=F\circ j$ and $\hat{f}|_U=\hat F\circ j$.

 One may expect the nowhere flat submanifolds
$f\colon\,M^n\to\R^{n+2}$ whose second fundamental form is ``as simple as possible", namely,  with constant
relative nullity index $\nu_f=n-2$, to be easily deformable. However, we believe that, as part of a general
trend, submanifolds in this class are generically rigid as happens in the hypersurface situation.  These rank
two  submanifolds where divided into \cite{DF3} into three classes: elliptic, parabolic and hyperbolic.  It was
shown there that  the elliptic and the nonruled parabolic ones are genuinely rigid.  The  ruled parabolic
submanifolds admit isometric immersions as hypersurfaces and have many isometric deformations. Examples of
hyperbolic submanifolds that are not genuinely rigid where discussed in \cite{DF3} but there are not general
results for this class yet.

In this paper, we generalize \cite{DF3} by showing that the nonruled parabolic submanifolds are not only genuinely
rigid but are, in fact,  isometrically rigid. In particular, this unexpected result provides the first known
examples of locally rigid submanifolds of rank and codimension two. Observe that by our result  any nonruled
parabolic submanifold cannot be locally realized as a hypersurface of a Sbrana-Cartan hypersurface. This is
certainly not the case for the elliptic submanifolds. The remaining of the paper is devoted to a local parametric
classification of all parabolic submanifolds.

\section{Preliminaries}

Let $f\colon\,M^n\to\R^{n+2}$ denote an isometric immersion with codimension two into Euclidean space
of a Riemannian manifold of dimension $n\geq 3$. We denote its second fundamental with values in the normal bundle by
$$
\alpha_f\colon\,TM\times TM\to T_f^\perp M.
$$
The shape operator $A^f_\xi\colon\, TM\to TM$ for any
$\xi\in T_f^\perp M$ is defined by
$$
\<A^f_\xi X,Y\>=\<\alpha_f(X,Y),\xi\>.
$$

We assume throughout the paper that $f$ has constant rank $2$. This condition is denoted by $rank_f=2$ and means that the relative nullity subspaces $\Delta(x)\subset T_xM$
defined by
$$
\Delta(x)=\{X\in T_xM : \alpha_f(X,Y)=0 \,;\, Y\in T_xM\}
$$
form a tangent subbundle of codimension two. Equivalently, the index of relative nullity $\nu_f(x)=\dim\Delta(x)$ satisfies $\nu_f(x)=n-2$.
It is a standard fact  that the relative nullity distribution is integrable and that the $(n-2)$-dimensional leaves are totally geodesic submanifolds of the manifold and the ambient space.\vspace{1ex}

The following fact was proved in \cite{DF3}.

\begin{proposition}\label{posto}\po
Let $f\colon\,M^n \rightarrow \R^{n+2}$ be an isometric immersion  of $\rank_f=2$. Assume that $M^n$ has no open flat subset. Then, given an isometric immersion $g\colon\,M^n \rightarrow \R^{n+2}$, there exists an open dense subset of $M^n$ such that along each connected component $V$ we have:
\begin{itemize}
\item[(i)] $\rank_g=2$ and $\Delta_g=\Delta_f$, or
\item[(ii)] $\rank_g=3$  and $g|_V=k\circ h\colon\,V\to\R^{n+2}$ is a composition of isometric immersions  $h\colon\,V\to U$ and $k\colon\,U\to\R^{n+2}$ where $U\subset\R^{n+1}$ is open.
\end{itemize}
\end{proposition}

The simplest submanifolds $f\colon\,M^n \rightarrow \R^{n+2}$ of $\rank_f=2$ are  called \textit{surface-like}. This means that
$f(M)\subset L^2 \times \R^{n-2}$ with $L^2 \subset\R^4$ or $f(M)\subset CL^2\times \R^{n-3}$ where $CL^2 \subset \R^5$ is a cone over a spherical surface $L^2\subset\mathbb{S}^4$.
The  following characterization in terms of the splitting tensor is well-known; see \cite{dft} or \cite{DG2}.
Recall that associated to the relative nullity foliation the
\textit{splitting tensor} $C$ is defined as follows: to each vector $T\in\Delta$ corresponds the endomorphism  $C_T$ of $\Delta^\perp$
given~by
$$
C_T X=-\left(\nabla_X T\right)_{\Delta^\perp}.
$$

\begin{proposition}\po\label{CT}
Let $f\colon\,M^n\rightarrow\R^{n+2}$ be an isometric immersion of  $\rank_f=2$. Assume that $C_T=\mu(T)I$ for any $T \in \Delta$. Then
each point has a neighborhood where $f$ is surface-like.
\end{proposition}

\section{Parabolic submanifolds}

A submanifold $f\colon\,M^n\rightarrow \R^{n+2}$ of $\rank_f =2$ is called  \textit{parabolic} if we have:
\begin{itemize}
\item [(i)] At any $x\in M^n$, we have
$$
T_{f(x)}^\perp M=\spa\{\alpha_f(X,Y)\,;\,X,Y\in T_xM\}.
$$
\item [(ii)] There is a nonsingular (asymptotic) vector field
$Z\in\Delta^\perp$ such that
$$
\alpha_f (Z,Z)=0.
$$
\end{itemize}

Clearly, if $f$ is parabolic there is
a smooth orthonormal frame $\{\eta_1,\eta_2\}$ of $T_f^\perp M$ such that
the shape operators have the form
\be\label{forma}
A^f_{\eta_1}|_{\Delta^\perp}=\left[\begin{array}{cc}
  a & b \\
  b & 0
\end{array}\right]\quad\textrm{and}\quad A^f_{\eta_2}|_{\Delta^\perp}=\left[\begin{array}{cc}
  c & 0 \\
  0 & 0
\end{array}\right]
\ee
with $b,c\in C^\infty(M)$ nowhere vanishing. In particular,
the asymptotic vector field $Z$ is unique up to sign.

\newpage

It is well-known (cf.\ \cite{DG2} or \cite{dft}) that the
differential equation
\be\label{dif}
\nabla _T A^f_\xi|_{\Delta^\perp}=A^f_\xi|_{\Delta^\perp} \circ C_T
+ A_{\nabla^\perp_T\xi}|_{\Delta^\perp}
\ee
holds for any $T\in\Delta$ and $\xi\in T_f^\perp M$. In particular, $A^f_\xi|_{\Delta^\perp} \circ C_T$ is symmetric. From this and
$(ii)$ it follows easily that
\be\label{r}
C_T=\left[\begin{array}{cc}
m & 0 \\
n & m
\end{array}\right]
\ee
for any  $T\in\Delta$. From Proposition \ref{CT}, we conclude
that $f$ is surface-like if and only  if $n=\<C_TX,Z\>=0$ for any $T\in\Delta$.\vspace{1ex}

Recall that an isometric immersion   $f\colon\,M^n\rightarrow\R^{n+2}$ is ruled when there is a foliation by open subsets of  $(n-1)$--dimensional affine subspaces of $\R^{n+2}$. Observe that  a ruled submanifold in codimension $2$ does not have to be parabolic. In fact, generically we have $\rank_f=3$.

Ruled parabolic submanifolds are never locally  isometrically rigid as seen in the converse part of the following result in \cite{DF3}. In fact, the direct statement will be an important element in the proof of our main result.

\begin{proposition}\po\label{ruled} If $f\colon\,M^n\rightarrow\R^{n+2}$ is a  ruled parabolic submanifold. If $M^n$ is simply connected then it admits an isometric immersion as a ruled hypersurface in $\R^{n+1}$.

Conversely, let $g\colon\,M^n\rightarrow\R^{n+1}$ be a simply connected ruled hypersurface without flat points. Then, the family of ruled parabolic submanifolds $f\colon\,M^n\rightarrow\R^{n+2}$ is parametrized by the set of ternary smooth arbitrary functions in an interval.
\end{proposition}

\section{The main result}

The following result generalizes the one in \cite{DF3} and provides
the first known examples of locally rigid submanifolds of rank and codimension two.

\begin{theorem}\po\label{main} Let $f\colon\,M^n \rightarrow \R^{n+2}$, $n\ge 3$, be a parabolic submanifold  neither ruled nor surface-like on any open subset of $M^n$. Then $f$ is  isometrically rigid.
\end{theorem}

 The following is a key ingredient of the proof of the theorem.

\begin{proposition}\po\label{teo:hiper}
Let $f\colon\,M^n\rightarrow \R^{n+2}$ be a simply connected  parabolic submanifold. If $f$ is not of surface type in any open subset of  $M^n$ and admits an isometric immersion as a hypersurface of $\R^{n+1}$ then $f$ is ruled.
\end{proposition}

\proof Let $g\colon\,M^n\rightarrow \R^{n+1}$ be an isometric immersion.
We denote by $N$ its  Gauss map. Given $x\in M^n$, let  $\beta\colon\,T_xM\times T_xM\rightarrow\Les^2$  be the symmetric bilinear form
$$
\beta(Y,V)=\< A^f_{\eta_1}Y,V\>e_1+\< A^g_N\,Y,V\>e_2
$$
where $\eta_1$ is as in (\ref{forma}) and $\{e_1,e_2\}$ is an orthonormal frame for the Lorentzian plane $\Les^2$ such that $\|\e_1\|^2=1=-\|e_2\|^2$ and $\<\e_1,e_2\>=0$. Then $\beta$ is flat since, from the Gauss equations for $f$ and $g$,  we easily see  that
$$
\<\beta(X,Y),\beta(V,W)\>
-\<\beta(X,W),\beta(V,Y)\>=0.
$$

We claim that $\Delta_g(x)=\Delta_f(x)$.  Since parabolic submanifolds have no flat points it follows that
$\dim\,\Delta_g(x)\le n-2$. If $\Delta_g(x)\neq\Delta_f(x)$,
it is easy to see that
$$
S(\beta)=\spa\{\beta(Y,V)\;;\; Y, V\in T_xM\}
$$
satisfies $S(\beta)=\Les^2$.  From  Corollary 1 in \cite{Mo}, we obtain that
$$
N(\beta)=\{Y\in T_xM : \beta(Y,V)=0\;;\;V\in T_xM\}
$$
satisfies $\dim N(\beta)=n-2$. Since
$N(\beta)=\Delta_g\cap\Delta_f$, this is a contradiction and proves the claim.

Set
$$
A^g_N|_{\Delta^\perp}
=\left[\begin{array}{cc}
\bar{a} & \bar{b}\\
\bar{b} & \bar{c}
\end{array}\right].
$$
Using  (\ref{r}) we conclude from the symmetry of
$$
A^g_N \circ C_T=\left[\begin{array}{cc}
\bar{a}m+\bar{b}n &\bar{b}m \\
\bar{b}m+\bar{c}n & \bar{c}m
\end{array}\right]
$$
that $\bar{c}n=0$. It follows from Proposition \ref{ruled} and our assumption that $f$ is nowhere surface-like that $\bar{c}=0$.
In particular, we have from the Gauss equations for $f$ and $g$ that
we can choose an orientation for $g$ such  that $\bar{b}= b$.

Taking the $Z$-component of the Codazzi equations for $A^f_{\eta_1}$
and $A^g_N$ give
\be\label{first}
2b\<\nabla_X X,Z\>-a\< \nabla_Z X,Z\> -Z(b)=0
\ee
and
$$
2b\<\nabla_X X,Z\>-\bar{a}\< \nabla_Z X,Z\>-Z(b)=0.
$$
Thus,
$$
\label{as}
(a-\bar{a})\< \nabla_Z Z,X\>=0.
$$

Suppose  that $\<\nabla_Z Z,X\>=0$ in an open subset $U$ of $M^n$.  From (\ref{r})
we have
\be\label{later}
\<\nabla_ZT,X\>=-\<C_TZ,X\>=0.
\ee
Then taking the  $Z$-component of the Codazzi equations for $g$ applied to $Z,T$~gives
$$
\<\nabla_TZ,X\>=0.
$$
It follows from the above that distribution $\spa\{Z\}\oplus\Delta$ is totally geodesic on $U$. But then $f$ is ruled on $U$ as we wanted.

Assume now that $a=\bar{a}$ on an open subset. Taking the $X$-component of the same Codazzi equations as above gives
\be\label{second}
X(b)-a\<\nabla_XZ,X\>-Z(a)+2b\<\nabla_Z
X,Z\>+c\<\nabla^\perp_Z\eta_1,\eta_2\>=0
\ee
and
$$
X(b)-a\<\nabla_XZ,X\>-Z(a)+2b\<\nabla_Z X,Z\>=0.
$$
Thus $\<\nabla^\perp_Z\eta_1,\eta_2\>=0$.
The Codazzi equation for $A^f_{\eta_2}$ applied to $X,Z$ yields
\be\label{igual}
c\< \nabla_ZZ,X\>
=b\<\nabla^\perp_Z\eta_1,\eta_2\>.
\ee
Again $\<\nabla_ZZ,X\>=0$, and the proof follows.\qed \vspace{1,5ex}

\noindent{\it Proof of Theorem \ref{main}}. Let $g\colon\,M^n\rightarrow\R^{n+2}$ be an isometric immersion. Since $f$ is nowhere ruled, Proposition \ref{teo:hiper} asserts that there is no local isometric immersion of $M^n$ as a hypersurface in $\R^{n+1}$. From Proposition \ref{posto}, we obtain $\rank_g=2$ and that $\Delta_g=\Delta_f$.

Take  $\{\eta_1,\eta_2\}$ as in (\ref{forma}) and recall that (\ref{r}) holds for any $T\in\Delta$.  Since
$A^g_{\bar{\eta}}\circ C_T$ is symmetric for any $\bar{\eta}\in T_g^\perp M$, it follows as before that $Z$ is asymptotic for~$g$. Therefore, we may fix a orthonormal  base $\{\bar{\eta}_1,\bar{\eta}_2\}$ of $T_g^\perp M$ such that
$$
A^g_{\bar{\eta}_1}|_{\Delta^\perp} =\left[\begin{array}{cc}
 \bar{a}& b \\
  b & 0
\end{array}\right]\quad \mbox{and}\quad
A^g_{\bar{\eta}_2}|_{\Delta^\perp} =\left[\begin{array}{cc}
  \bar{c} & 0 \\
  0 & 0
\end{array}\right].
$$
Since (\ref{first}) holds for both immersions, we have
$$
(a-\bar{a})\< \nabla_ZZ, X\>=0,
$$
and conclude that  $a=\bar{a}$.  Similarly, we have from (\ref{second}) that
$$
c\<\nabla^\perp_Z\eta_1,\eta_2\>
=\bar{c}\<\nabla^\perp_Z\bar{\eta}_1, \bar{\eta}_2\>.
$$
Now computing the $Z$-component of the Codazzi equations for $A^f_{\eta_2}$ and
$A^g_{\bar{\eta}_2}$ gives
$$
\bar{c}\<\nabla^\perp_Z\eta_1,\eta_2\> =c\<\nabla^\perp_Z\bar{\eta}_1,\bar{\eta}_2\>.
$$
It follows that
$$
c=\bar{c} \;\;\;\;\mbox{and}\;\;\;\;
\<\nabla^\perp_Z\eta_1,\eta_2\> =\<\nabla^\perp_Z\bar{\eta}_1,\bar{\eta}_2\>.
$$
Then taking the $X$-components yields
$$
\<\nabla^\perp_X\eta_1,\eta_2\> =\<\nabla^\perp_X\bar{\eta}_1,\bar{\eta}_2\>.
$$
We conclude from the fundamental theorem of submanifolds (cf.\ \cite{da}) that $f$ and $g$ are congruent by a rigid motion of the ambient space.\qed

\section{Ruled parabolic surfaces}

In this  section, we give a parametric description of all ruled  submanifolds in $\R^{n+2}$ that are parabolic.
\vspace{1,5ex}

Let $v\colon\,I\to\R^{n+2}$ be a smooth curve parametrized by arc length in an interval $I$ in $\R$. Set  $e_1=dv/ds$ and let $e_2,\ldots, e_{n-1}$ be  orthonormal normal vector fields along $v=v(s)$ parallel in the normal connection of  $v$ in $\R^{n+2}$. Thus,
\be\label{condi0}
\frac{de_j}{ds}=b_je_1,\;\;\; 2\le j\le n-1,
\ee
where $b_j\in C^\infty(I)$. Set $\Delta=\spa\{e_2,\ldots, e_{n-1}\}$ and let $\Delta^\perp$ be the orthogonal  complement in the normal bundle of the curve.
Take a smooth unit vector field $e_0\in\Delta^\perp$ along $v$ such that
$$
P=\spa\{e_0,(de_1/ds)_{\Delta^\perp}\}\subset\Delta^{\perp}
$$ satisfies:
\be\label{dim} \dim P = 2 \ee and $P$ is nowhere parallel in $\Delta^\perp$ along $v$, that is, \be\label{dim2}
\spa\{(de_0/ds)_{\Delta^\perp},  (d^2e_1/ds^2)_{\Delta^\perp}\} \not\subset P. \ee

We parametrize a ruled submanifold $M^n$ by
\be\label{para} f(s,t_1,\ldots,t_{n-1})={c(s)}+\sum_{j=1}^{n-1}t_j e_j(s)
\ee
where $(t_1,\ldots,t_{n-1})\in\R^{n-1}$ and $c(s)$ satisfies $dc/ds=e_0$.
To see that $f$ is parabolic, first observe that
$$
TM=\spa\{f_s\}\oplus\spa\{e_1\}\oplus\Delta
$$
where
$f_s=e_0+t_1de_1/ds+\sum_{j\ge 2}t_jb_je_1$.
Consider the orthogonal decomposition
\be\label{eum}
\left(\frac{de_1}{ds}\right)_{\Delta^\perp}= a_1e_0 + \eta.
\ee
Thus $\eta(s)\neq 0$ for all $s\in I$ from (\ref{dim}). Hence,
\be\label{tangente}
TM=\spa\{e_0+t_1(a_1e_0+\eta)\}\oplus\spa\{e_1\} \oplus\Delta,
\ee
and it follows easily that $f$ is regular at any point .

Since $f_{st_j}=b_je_1\in TM,\; 2\le j\le n-1$,
we have that $\Delta\subset\Delta_f$. It follows easily from (\ref{eum}), (\ref{tangente}) and $\eta(s)\neq 0$ that
$$
f_{st_1}=\frac{de_1}{ds}\not\in TM.
$$
It is easy to see that $f_{ss}\not\in \spa\{f_{st_1}\}\oplus TM$, i.e.,  $\dim N_1^f=2$, is equivalent to
$$
\left(\frac{de_0}{ds}\right)_{\Delta^\perp} +t_1\left(\frac{d^2e_1}{ds^2}\right)_{\Delta^\perp} \not\in P.
$$
It follows easily that $\Delta=\Delta_f$ and that $f$ is parabolic in at least an open dense subset of $M^n$.

Let $f\colon M^n\to\R^{n+2}$ be a ruled parabolic submanifold and
$\{e_2,\ldots,e_{n-1}\}$ an orthonormal frame
for $\Delta_f$ along a integral curve $c=c(s),\; s\in I$, of the
unit vector field $X$ orthogonal to the rulings.  Without loss of generality (see Lemma~2.2 in \cite{BDJ}) we may assume that
$$
\frac{de_j}{ds}\perp\Delta_f,\;\;2\le j\le n-1.
$$
Now parametrize $f$ by (\ref{para}), where $e_0=X$ and $e_1=Z$.
That $f_{st_j}\in TM$ implies
\be\label{condi11}
\frac{de_j}{ds}\in\spa\{e_1,f_s\},\;\;\; 2\le j\le n-1.
\ee
Taking $t_1=0$, we obtain that
\be\label{condi22}
\frac{de_j}{ds}=a_je_0+b_je_1,\;\;\; 2\le j\le n-1,
\ee
where $a_j,b_j\in C^\infty(I)$.  Since  $\dim N^f_1=2$, we have
$$
\frac{de_1}{ds}= a_1e_0 + (de_1/ds)_{\Delta} + \eta
$$
where $\eta\perp \spa\{e_0,e_1\}\oplus\Delta$ satisfies
$\eta(s)\neq 0$. Thus (\ref{condi11}) reduces to
$$
a_je_0\in\spa\{(1+t_1a_1+\ldots
+t_{n-1}a_{n-1})e_0+t_1\eta\} ,\;\;\; 2\le j\le n-1.
$$
Therefore $a_j=0$. From (\ref{condi22}) we have
$de_j/ds=b_je_1$ for $2\le j\le n-1$.\vspace{1,5ex}

We have proved the following result.

\begin{proposition}\po Given a smooth curve $c\colon\, I\subset \R\rightarrow \R^{n+2}$ consider orthonormal fields $e_0=dc/ds, e_1(s),\ldots,e_{n-1}(s)$ satisfying (\ref{condi0}), (\ref{dim}) and (\ref{dim2}) at any point. Then, the ruled submanifold parametrized by
\be\label{paramregrada}
f(s,t_1,\ldots,t_{n-1})={c(s)}+\sum_{j=1}^{n-1}t_j e_j(s)
\ee
with $(t_1,\ldots,t_{n-1})\in\R^{n-1}$ is  parabolic in an open dense subset of $M^n$.
Conversely, any ruled parabolic submanifold can be parametrized as in (\ref{paramregrada}).
\end{proposition}

\section{Nonruled parabolic surfaces}

In this section, we provide a parametrically description of all nonruled Euclidean parabolic submanifolds. \vspace{1,5ex}

Let $L^2$ be a Riemannian manifold endowed with a  global coordinate system $(x,z)$. Let $g\colon\,L^2\to\R^N$,
$N \geq 4$, be a surface  whose coordinate functions are linearly independent  solutions of  the parabolic equation
\be\label{PDE}
\frac{\partial^2 u}{\partial z^2}+W(u)=0
\ee
where $W\in TL$. In terms of the Euclidean connection, we
have
$$
\tilde\nabla_Zg_*Z+g_*W=0
$$
where $Z=\partial/\partial z$.
Thus, the coordinate filed $Z$ is asymptotic, i.e., the second fundamental form  of $g$ satisfies
\be\label{ass}
\alpha_g(Z,Z)=0,
\ee
and also $W=-\nabla_ZZ$. In particular, the coordinate functions satisfy
\be\label{PDE1}
\hess_u(Z,Z)=0.
\ee
Conversely, if $f\colon\,L^2\to\R^N$ is a surface with a coordinate system $(x,z)$ such that $\partial/\partial z=Z$
satisfies (\ref{ass}),  then all coordinate functions of  $f$ satisfy (\ref{PDE}) with $W=-\nabla_ZZ$.

Given a surface $g\colon\,L^2\to\R^N$, $N \geq 4$, we denote $$
N^g_1(x)=\spa\{\alpha_g(X,Y)\,;\,X,Y\in T_xL\}.
$$
We call  $g$ a \textit{parabolic} surface if the following conditions hold.
\begin{itemize}
\item [(i)]  $\dim N^g_1(x)=2$ at any $x\in L^2$.
\item [(ii)] There is a nonsingular vector field
$Z\in TL$ such that
$\alpha_g (Z,Z)=0.$
\end{itemize}

Let
$h\colon L^2\to \R^{N}$ be a smooth map satisfying
\be\label{section}
h_*(TL)\subset T_g^\perp L.
\ee
Set
$h=U+\delta$
where $U\in TL$ and $\delta \in T_g^{\perp}L$. Given $Y\in TL$,
we have
$$
h_\ast(Y)=\nabla_YU-A^g_\delta(Y)
+\alpha_g(Y,U)+\nabla^{\perp}_Y\delta.
$$
It follows that (\ref{section}) is equivalent to
\be\label{pde0}
\nabla_Y U=A^g_\delta Y,\;\;\mbox{for any}\;\; Y\in TL.
\ee
In particular, the map $(Y,X)\mapsto\<\nabla_Y U, X\> $ is symmetric, i.e., the one-form $U^*$ is closed. Thus
$U=\nabla\varphi$ for some function $\varphi\in C^\infty(L)$. We obtain from (\ref{pde0}) that
$$
\hess_\varphi=A^g_\delta
$$
and hence $\varphi$ satisfies (\ref{PDE1}).
Let $\Lambda$ denote
the orthogonal  complement of $N_1^g$ in $T_g^\perp L$, that is,
$$
T_g^\perp L=N_1^g\oplus \Lambda.
$$

We have just proved the direct statement in the following result since the proof of the converse is immediate.

\begin{proposition}\po Let $g\colon\,L^2\to\R^{N}$ be a parabolic surface. Then any smooth map $h\colon\, L^2\to\R^N$ satisfying (\ref{section}) has the form
\be\label{secao}
h_\varphi= g_*\nabla\varphi+\gamma_1+\gamma_0
\ee
where the function $\varphi$ satisfies (\ref{PDE})(or (\ref{PDE1})), the section $\gamma_1 \in N_1^g$ is unique such that $A^g_{\gamma_1}=\hess_\varphi$ and $\gamma_0$ is any section of $\Lambda$. Conversely, any function $h_\varphi$ as in  (\ref{secao}) satisfies (\ref{section}).
\end{proposition}

Let $f\colon\,M^n\to\R^{n+2}$ be parabolic. For simplicity, we assume that $f$ does not split a Euclidean factor.
Since we are working locally, we assume that $M^n$ is the saturation of a transversal section $L^2$ to $\Delta$. In relation to the following definition recall that the normal bundle $T_f^\perp M$ is constant along $\Delta$ in $\R^{n+2}$.

\begin{definition}\po{ \em
Let $f\colon\,M^n\to\R^{n+2}$ be a parabolic submanifold and let $L^2$ be as above. We call a \textit{polar surface} associated  to $f$ any immersion  $g\colon\,L^2\to\R^{n+2}$ satisfying $T_{g(x)}L=T_{f(x)}^\perp M$ up to parallel identification in $\R^{n+2}$.}
\end{definition}

\begin{proposition}\po
Any parabolic submanifold $f$ admits locally a polar surface.
Moreover, any polar surface $g$ associated to $f$ is parabolic and  substantial.
\end{proposition}

\proof Let $\eta_1, \eta_2 \in T_f^\perp M$ be orthonormal vector fields such that (\ref{forma}) holds. We claim that
they are parallel along $\Delta$ in the normal connection.
It follows from (\ref{dif}) that
$$
\nabla _T A^f_\xi=A^f_\xi \circ C_T
$$
is satisfied if $\xi\in T_f^\perp M$ is parallel along $\Delta$. Given $x\in M^n$, let $\gamma$ be a geodesic with
$\gamma(0)=x$ contained  in the leaf of relative nullity tangent to $\Delta(x)$.
If $\delta_t$ is the parallel transport of $\eta_x$ along $\gamma$, we have
$$
\nabla_{\gamma'}A^f_{\delta_t}=A^f_{\delta_t}\circ C_{\gamma'}.
$$
Hence,
$A^f_{\delta_t}=A^f_{\eta_x}e^{\,\int_0^t C_{\gamma'}d\tau}$.
Thus  $A^f_{\delta_t}$ has constant rank along $\gamma$.
It follows that $A_\xi$ has constant rank along the leaves of relative nullity.
Now observe that $\eta_1$ is the unique (up to sign) unit vector field in $T_f^\perp M$ such that the corresponding shape operator $A^f_{\eta_1}$ has rank one.  Thus $\eta_1$ and hence $\eta_2$ are parallel as claimed.

Next, we claim that also the orthonormal frame $\{X,Z\}$ in $\Delta^{\perp}$ as in (\ref{forma}) is parallel along $\Delta$. In fact, using (\ref{later}) and the Codazzi equation we obtain
$$
\nabla_T^\perp\alpha_f(Z,X) =\<\nabla_TZ,X\>\alpha_f(X,X)-\<\nabla_Z T,Z\>\alpha_f(Z,X).
$$
Since $\eta_1$ is colinear with $\alpha_f(Z,X)$ and parallel along $\Delta$, the claim follows easily.

Let $U,V\in TL$ be such that
$$
Z=U+\delta_1,\;\;\; X=V+\delta_2
$$
where $\delta_1,\delta_2\in\Delta$. Now let $(u,v)$ be a coordinate system in $L^2$ such that
$$
\partial u=:\partial/\partial u=\lambda_1U,\;\;\;
\partial v=:\partial/\partial v=\lambda_2V
$$
where $\lambda_1,\lambda_2\in C^\infty(L)$.
We show next that there exist linearly independent one-forms
$\theta_1, \theta_2$  such that the differential equation
\be\label{eqdf}
dg=\theta_1 \eta_1+\theta_2\eta_2
\ee
is integrable.  Set
\be\label{eqdf2}
\theta_1=c\lambda_2\phi dv \quad\mbox{and}
\quad \theta_2=b\lambda_1\phi du+\sigma dv
\ee
where $\phi,\sigma \in C^\infty(L)$ and $b,c$ denote the restriction to $L^2$ of the functions in $M^n$ defined in (\ref{forma}). The  integrability
condition of  (\ref{eqdf}) is
\bea 0\!\!&=&\!\! d\theta_1\eta_1+d\theta_2\eta_2
+\theta_1 \land d\eta_1+\theta_2\land d\eta_2\\
\!\!&=&\!\!d\theta_1\eta_1+d\theta_2\eta_2
-\left(c\lambda_2\phi\partial \eta_1/\partial u
-b\lambda_1\phi\partial\eta_2/\partial v
+\sigma\partial\eta_2/\partial u\right)dV
\eea
where $dV$ is the volume element of $L^2$.
We have
$$
\tilde\nabla_{\partial u}\eta_j
=\tilde\nabla_{\lambda_1(Z-\delta_1)}\eta_j
=\lambda_1\tilde\nabla_{Z}\eta_j
\;\;\;\mbox{and}\;\;\;
\tilde\nabla_{\partial v}\eta_j
=\tilde\nabla_{\lambda_2(X-\delta_2)}\eta_j
=\lambda_2\tilde\nabla_{X}\eta_j.
$$
Thus, we obtain from (\ref{forma}) that
$$
(c\lambda_2\phi\partial \eta_1/\partial u
-b\lambda_1\phi\partial\eta_2/\partial v
+\sigma\partial\eta_2/\partial u)_{TM}=0.
$$
Then,
$$
c\lambda_2\phi\partial \eta_1/\partial u
-b\lambda_1\phi\partial\eta_2/\partial v
+\sigma\partial\eta_2/\partial u
=e\eta_1+\ell\eta_2
$$
where $e,\ell\in C^\infty(L)$.
To conclude, we observe that the integrability condition  now follows from the integrability of the system
$$
\left\{\begin{array}{lll}
(c\lambda_2\phi)_u \!\!&=&\!\!e\vspace*{1ex} \\
\sigma_u-(b\lambda_1\phi)_v \!\!&=&\!\!\ell.
\end{array}\right.
$$

It remains to see that $g$ is parabolic. It is clear that $N_1^g=\Delta^\perp$, and hence $\dim N_1^g=2$. It
follows from (\ref{eqdf2}) that $g_*(\partial u)=b\lambda_1\phi \eta_2$. Then, we obtain from (\ref{forma}) that
$\tilde{\nabla}_{\partial u}g_*(\partial u)\in TL$ as we wished.\qed

\begin{theorem}\po\label{theorem:polar} Let $g\colon\,L^2\to \R^{n+2}$ be a parabolic surface and let \mbox{$\Psi\colon\,\Lambda\to\R^{n+2}$} be the map defined by \be\label{paramet1}
\Psi(\delta)=h(x)+ \delta ,\;\; \delta \in \Lambda(x),
\ee
where $h\colon L^2\to\R^{n+2}$ satisfies (\ref{section}). Then,   $M^n=\Psi(\Lambda)$ is, at regular
points, a  parabolic submanifold
with polar surface $g$. Moreover, $\Psi$ is nonruled if and only
if $g$ is nonruled.

$\!\!\!$Conversely, any parabolic submanifold $f\colon\,M^n\to\R^{n+2}$  admits  a parametrization
(\ref{paramet1}) locally where $g$ is a polar surface of $f$.
\end{theorem}

\proof Since $g$ is parabolic there is an orthonormal
tangent frame $\{X,Z\}$ of $TL$ and an orthonormal normal
frame of  $\{\eta_1,\eta_2\}$ of $N_1^g$ such that
$$
A^f_{\eta_1}=\left[\begin{array}{cc}
  a & b \\
  b & 0
\end{array}\right]\quad\textrm{and}\quad A^f_{\eta_2}=\left[\begin{array}{cc}
  c & 0 \\
  0 & 0
\end{array}\right]
$$
with $b,c\in C^\infty(M)$ nowhere vanishing.

It follows from the symmetry of the tensor (see \cite{Sp})
$$
\alpha^2(V,Y,W)
=(\nabla_V^\perp\alpha_g(Y,W))_{(N_1^g)^\perp}
$$
that \be\label{l8} (\nabla_Z^\perp\alpha_g(X,Z))_{(N_1^g)^\perp}=0, \ee that is, $$
(\nabla_Z^\perp\eta_1)_{(N_1^g)^\perp}=0. $$

 Next we show that \be\label{l7} (\tilde\nabla_Zh_*(Z))_{TL}=0. \ee In fact, we have
$$
\<\tilde{\nabla}_Z h_*(Z),W\>
=\<\tilde{\nabla}_Z h,\alpha_g(W,Z)\>= \<\tilde{\nabla}_W \tilde{\nabla}_Z
h,Z\> =\<\tilde{\nabla}_Z \tilde{\nabla}_W h,Z\>=0
$$
since $Z$ is asymptotic.

Since $h$ satisfies (\ref{section}) it follows easily from the regularity assumption that $T_{\Psi(\delta_x)}M=T_g^\perp L(x)$. Moreover, we have that
$\Delta_{\Psi(\delta_x)}=\Lambda(x)$. Let us see that $\Psi$ is parabolic. Since
$$
\Psi_*(Z)=h_*(Z)+\tilde{\nabla}_Z\delta\in
T_{g(x)}^\perp L=T_{\Psi(\delta_x)}M,
$$
we have using (\ref{l8}) and (\ref{l7}) that
\begin{eqnarray}
\label{asymptotic} \<\tilde{\nabla}_Z\Psi_*(Z),W\>
\!\!&=&\!\!\<\tilde{\nabla}_Z h_*(Z)+\tilde{\nabla}_Z\delta,\alpha_g(W,Z)\>\\
\!\!&=&\!\!\<\tilde{\nabla}_Z h_*(Z),\alpha_g(W,Z)\>
+\<\delta,\nabla^\perp_Z\alpha_g(W,Z)\>\nonumber\\
\!\!&=&\!\! 0\nonumber.
\end{eqnarray}
Since $\eta_1\in T^\perp_g L = T_\Psi M$, and because $(\tilde{\nabla}_X \eta_1)_{TL}, (\tilde{\nabla}_Z
\eta_1)_{TL}$ are linearly independent and belong to $T^\perp_\Psi M$, then $\dim N_1^\Psi=2$.  Thus $\Psi$ is
parabolic.

 Let us see that $\Psi$ is not ruled. We have from (\ref{asymptotic}) that the asymptotic direction of $\Psi$ is
normal to $\eta_1$, thus colinear with $\eta_2$. Since  $\Psi$ is ruled if and only if $\<\tilde{\nabla}_Z\eta_2, \eta_1\>=0$, we have from (\ref{igual}) that $g$ is ruled.

For the converse, consider a polar surface $g\colon\,L^2\to \R^{n+2}$ of $f$. Since $f$ has no Euclidean factor
then $g$ is substantial. It is now easy to conclude that  $\Delta_f=\Lambda$ and $TM=T_g^\perp L$ along $L^2$.
Thus, $h=f_{\mid_{L}}$ satisfies (\ref{section}).\qed

\vspace{.5in} {\renewcommand{\baselinestretch}{1}
\hspace*{-20ex}\begin{tabbing} \indent\= IMPA -- Estrada Dona Castorina, 110
\indent\indent\=  Universidade da Beira Interior\\
\> 22460-320 -- Rio de Janeiro -- Brazil  \>
6201-001 -- Covilhã -- Portugal \\
\> E-mail: marcos@impa.br \> E-mail: pmorais@mat.ubi.pt
\end{tabbing}}

\end{document}